# Stochastic-based Optimal Daily Energy Management of Microgrids in Distribution Systems

Farhad Samadi Gazijahani
Electrical Engineering Department
Azarbaijan Shahid Madani University
Tabriz, Iran
f.samadi@azaruniv.ac.ir

Javad Salehi
Electrical Engineering Department
Azarbaijan Shahid Madani University
Tabriz, Iran
j.salehi@azaruniv.ac.ir

*Abstract*—Microgrid (MG) with different technologies in distributed generations (DG) and different control facilities require proper management and scheduling strategies. In these strategies, in order to reach the optimal management, the stochastic nature of some decision variables should be considered. Therefore, we proposes a new stochastic-based method for optimal daily energy management (SDEM) of MGs considering economic and reliability aspects. The optimization aim is to minimize overall operating cost, power losses cost, pollutants emission cost and cost of energy not supply (ENS). The network is assumed to be supplied by renewable and dispatchable generators and energy storage systems (ESS). The system uncertainties are considered using a set of scenarios and a scenario reduction method is applied to enhance a tradeoff between the accuracy of the solution and the computational burden. Cuckoo optimization algorithm (COA) is applied to minimize the objective function as an optimization algorithm. The effectiveness and efficiency of the proposed method are validated through extensive numerical tests on PG&E 69-bus test distribution system. The results show that the proposed framework can be considered as an efficient tool in optimal daily energy management of smart distribution networks.

## Nomenclature

| | |
|---|---|
| $s$ | Set of scenarios |
| $b$ | Branch |
| $P_{WT(v)}$ | Generated power at wind turbine |
| $P_{r,WT}$ | Rated power of WT |
| $P_D$ | Load demand |
| $v$ | Wind speed |
| $P_{PV}$ | Output power of the module at irradiance |
| $G_{ING}$ | Incident irradiance |
| $G_{STC}$ | Irradiance at standard test condition |
| $P_{STC}$ | Rated power of photovoltaic |
| $Z$ | Objective function |
| $T_c, T_r$ | Cell reference temperatures respectively |
| $k$ | Maximum power correction temperature |
| $\alpha, \beta$ | Shape and scale parameter of Beta |
| $\Gamma$ | Parameters of the Weibull model |
| $\Theta, \varrho$ and $\gamma$ | CHPs coefficients |
| $\lambda$ | Failure rate |
| $SOC$ | State of charge |
| $\sigma$ | Standard deviation |
| $\pi$ | Mean value of distribution function |
| $\zeta, \delta$ | Generator emission characteristics |
| $ENS$ | Energy not supply |
| $\eta$ | Efficiency |
| $\mu$ | Membership function |
| $C_{int}$ | Price of energy not supply |
| $C_{Ploss}$ | Price of losses |
| $N_{res}$ | Number of nodes isolated during fault |
| $N_{rep}$ | Number of nodes isolated during fault |
| $P_{res}$ | Loads are restored during fault |
| $P_{rep}$ | Loads are not restored during fault |
| $T_{res}$ | Duration of fault location and switching |
| $T_{rep}$ | Duration of the fault repair |
| $R_b$ | Resistance of the line b |
| $L_b$ | Length of the line b |
| $H$ | Weight Coefficients |

## I. Introduction

Nowadays, implementation of the MG concept because of its effects on improving network security and environmental aspects is expanding in the distribution network [1]. From the perspective of MG owners, obviously, economic operation of the MG is important. Since, MGs can participate in power markets and also provide some ancillary services, proper scheduling of the MG is essential from the main grid point of view [2]-[5]. Therefore a suitable strategy should be pursued for the MG operation. The MG consists of distributed generators as well as renewable energy sources (RESs). There is important difference between the optimization of a MG and conventional economic dispatch problem [6]. At present, there are some researches on MG optimal energy dispatch. In the grid-connected mode, the optimal energy dispatch is analyzed in MG [7] and the effect of time-of-use electricity price and electric energy transaction are considered. To consider the decentralized optimal power dispatch strategies in [8], the sharing of marginal cost of each power source is done through iterative and communication. In [9], the real time energy optimization scheduling method is proposed in independent operation mode of MG. The economic dispatch optimization problem is solved with different methods in literature [10].

For example there are different approaches to tackle this kind of optimization problem, such as finding multiple Pareto optimal solutions [11] or optimizing a single multi objective function obtained by a weighted sum of the targets [12]. The most used optimal power flow (OPF) solution techniques are the Lagrange–Newton method [13] and sequential quadratic programming (SQP) [14]. In the last few years, stochastic optimization methods have been applied to the OPF problem. Reference [15] presents an approach to solve a single objective OPF problem (to minimize the total operating cost in a power system) by means of particle swarm optimization algorithm. In [16], the multi objective economic-emission OPF problem is solved using a differential evolution algorithm.

It should be mentioned that the duty of each MG is to create balance between supply and demand. In other word, in simple economic dispatch problem power generation sources such as small-scale energy resources (SSERs) produce power in order to supply load demand and with increasing or decreasing fuel in diesel generators, the balance between supply and demand is provided. But considering MG concept in smart grids, the balance between generation and load is done through power exchanging between MGs as well as main grid so that the total cost of power generation in each MG as well as the total cost of power exchanging between MGs and main grid be optimized.

This paper proposes a new approach for probabilistic optimal energy management considering uncertainties in load and probabilistic modeling of generated power by renewable SSERs. In order to deal with these uncertainties, PDF is considered for power generation parameters. The problem is solved with cuckoo optimization algorithm. The main contributions of this paper to the research field are related to stochastic optimum energy management with considering reliability and economical aspects for MGs in distribution system that will benefit utilities, DG owners and electricity consumers.

The rest of this paper is organized as follows: Section 2 discusses the probabilistic analysis method and concepts. In Section 3, the generation and load models are explained. The problem is formulated in Section 4 and Section 5 explains the solution algorithms and Section VI discusses implementation and lastly, the conclusions are drawn in Section VIII.

## II. Probabilistic Analysis Method

At present, high penetration of SSERs into distribution grids affect operation and planning of the power systems. In WTs and PV power generation, wind speed, and solar radiation are prime energy sources, respectively. Because of stochastic behavior of wind speed and sun irradiance, power generation of the above energy resources undergoing significant uncertainties. Uncertainties analysis of the impact of SSERs such as wind turbine (WT) and photovoltaic (PV) units on current distribution systems based on deterministic methods is complicated scenario. The probabilistic analysis of distribution systems at presence of uncertainties is a powerful tool for optimal operation and of distribution systems. In probabilistic analysis, input data have probability distribution function (PDF) and these data can be described by cumulative distribution function (CDF). Consequently, the obtained results from probabilistic analysis are also presented in PDF and CDF forms. The uncertainty modeling techniques are applied in [17] that authors are performed a good review on uncertainty modeling approaches. In this paper, PDF is considered for input data such as load, generated power by SSERs, cost of transaction, and operation costs. Hence, the output results are represented in framework of PDF and CDF.

## III. System Model

### A. Modeling of DG Units

The model describing each DG is described as follows.

*a) WT*: The power generated by WT as a function of wind speed can be calculated by:

$$P_{WT} = \begin{cases} 0 & 0<V<P_{STG} \\ (A.V^2+B.V+C)*P_{rate} & V_{ci} \leq V \leq V_r \\ P_{rate} & V_r \leq V \leq V_{co} \\ 0 & V_{co} \leq V \leq \infty \end{cases} \quad (1)$$

The relationship between the generated power and wind speed of the WT can be demonstrated as Fig.1.

*b) PV*: The power generated by PV depends on the irradiance and the ambient temperature obtained by:

$$f(s) = \begin{cases} \frac{\Gamma(\alpha+\beta)}{\Gamma(\alpha)+\Gamma(\beta)} \times s^{(\alpha-1)}.(1-s)^{\beta-1} & , 0 \leq s \leq 1, \alpha \geq 0, \beta \geq 0 \\ 0 & , otherwise \end{cases} \quad (2)$$

$$\beta = (1-\mu).(\frac{\mu.(1+\mu)}{\sigma^2}-1) \quad (3)$$

$$\alpha = \frac{\mu.\beta}{1-\mu} \quad (4)$$

$$P_{pv} = P_{STG} * \frac{G_{ING}}{G_{STG}} * (1+k(T_C-T_{ref})) \quad (5)$$

*c) CHP*: The output power generated by CHP is controlled by an installed governor. The rate of fuel consumed by CHP can be expressed by a second-order polynomial in terms of the power generated.

$$L_{CHP}(P_{CHP}) = \Theta.P_{CHP}^2 + \varrho.P_{CHP} + \gamma \quad (6)$$

*d) ESS*: Energy storage system enhances flexibility in power generation, reliability, and consumption. The state of charge (SOC) value at time $t$ is related to SOC value at time $t$-1 by the following equation.

$$SOC_t = SOC_{t-1} + \eta_{ch} P_{ch} \Delta t - \frac{1}{\eta_{dis}} P_{dis} \Delta t \quad (7)$$

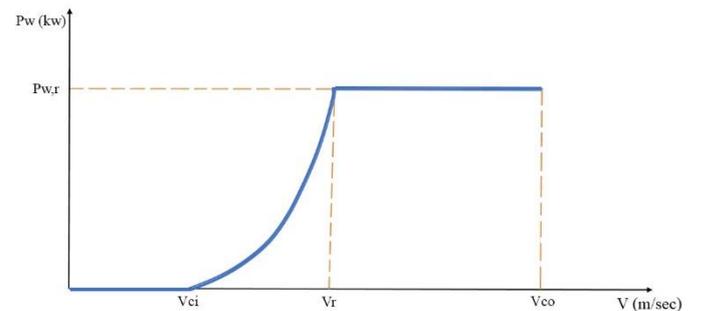

Fig. 1. Output power of wind generation unit related to wind speed.

## B. Load Demand Model

Because of variable nature of load, it is suitable to model the stochastic behavior of the load. This modeling can be achieved either by measured data analysis or by application of mathematical model. The MG load behavior is modeled as a normal distribution function. In this case, the MG load is defined by a mean value ($\mu$) and standard deviation ($\sigma$) as:

$$f(p_l) = \frac{1}{\sqrt{2\pi} \times \sigma} exp(-\frac{(p_l - \mu)^2}{2 \times \sigma^2}) \quad (8)$$

## IV. Proposed Model

### A. Cost evaluation

The economic-environment energy management between resources is obtained by minimizing $F1$. This function, simultaneously, minimizes fuel cost, operation cost, pollutants emission cost and losses cost. Renewable energy resources, such as WT and PV, do not incur expenses associated with pollution and fuel. Accordingly, their maintenance cost is considered a fixed value. Such value, therefore, is independent of the power generated. In the grid connected mode, WT and PV operate at their maximum power point and the rest of the power demand is provided by other distributed generation (DG) can be calculated by $F1$. The optimal values, obtained by minimizing $F1$, provide the amount of power generated by each DG, as well as power sold to or purchased from the grid on an hourly basis.

$$F1 = \sum_{t=1}^{T} \sum_{s=1}^{S} \sum_{n=1}^{N} \left[ f_{n,t}^s + OM_{n,t}^s + E_{n,t}^s + PLC_{loss} \right] \quad (9)$$

*a) Fuel Consumption Function*: Fuel expenses can be calculated as a function of power generated by the $n$th DG unit. The second term in [18], is an income coming from selling the thermal energy of CHPs. CHPs thermal energy is used to supply the thermal load. In the case of shortage of the thermal energy, the boiler supplies the rest of the required thermal energy.

$$f_{CHP}(P_{CHP}) = (C_{gasCHP} * \frac{P_{CHP}}{\mu_{CHP}}) - (C_{th}(\Theta_{H/E})_{CHP} * P_{CHP}) \quad (10)$$

*b) O&M Cost Function:* The operation and maintenance (O&M) cost of the $n$th unit as a function of the power generated can be obtained by:

$$OM_n(p_n) = k_{OM} * p_n \quad (11)$$

*c) Pollution Function:* The main source of air pollution and greenhouse gas emissions are seen to be fossil fuel consumption producing gases such as $CO_2$, $SO_2$, $NO$, and $CO$. The emission of each gas is converted to a corresponding expense through multiplication by a coefficient as indicated by:

$$E(P_G) = \sum_{i=1}^{Ne} (\alpha_i + \beta_i P_{G_i} + \gamma_i P_{G_i}^2) + \zeta_i exp(\lambda_i P_{G_i}) \quad (12)$$

*d) Power losses cost:* Calculation of power losses implies solving the load flow problem. In this paper backward-forward approach is applied for calculation of power losses.

$$PLC(t) = ((\sum_{b=1}^{Nb} [(P_b^2 + Q_b^2) \times R_b / V_b^2]) * C_{P_{loss}}) \quad (13)$$

### B. Reliability evaluation

Reliability evaluation of MGs in distribution systems with renewable energy resources (RERs) has attracted the interest of researchers worldwide. The reliability assessment process is complicated by different output characteristics of DGs in MGs. In this paper, the reliability evaluation is done for MGs under probabilistic and uncertain behaviors of MG components. One of the reliability indices is energy not supplied (ENS), in which it describes the capacity of loss of load (kWh). ENS is calculated based on power availability of each bus. Let $ENS_i$ be the loss of load obtained for the $i$th contingency, with a probability of $prob_i$. Then the average energy not supply or loss of load expectation (AENS or LOLE) is given by the following:

$$AENS = \sum_{I}^{N_c} ENS \times prob_i \quad (14)$$

The reliability of the network or energy index of reliability (EIR) is then given by the following:

$$EIR = 1 - \frac{AENS}{P_I} \quad (15)$$

$$C_{AENS} = \sum_{s=1}^{S} \sum_{b=1}^{Nb} Cint_b^s \lambda_b^s L_b (\sum_{res=1}^{N_{res}} P_{res}^s T_{res}^s + \sum_{rep=1}^{N_{rep}} P_{rep}^s T_{rep}^s) \quad (16)$$

In addition to mentioned indices, two other metrics of reliability which are proposed in [19], are calculated as Renewable Energy Penetration (REP) is one of metrics that describes the percentage of demand covered by renewable energy (WT and PV units) in a MG. All indices are calculated for MGs. In this paper, we set the MGs cost of interruption in 1.5 USD/kWh as it is described for household consumption. We use the following equation to calculate the total interruption cost in MGs for the reliability index on a consumer side:

$$IC_{iday} = H_c \times C_{AENS} \quad (17)$$

$$F2 = \sum_{MG} (IC_{dey,MG}) \quad (18)$$

### C. Objective Function

In conclusion, economic and reliability viewpoints which have been explained in previous parts are considered in two unique objective function that formulated below:

$$Z = [H_1 * F1 + H_2 * F2] \quad (19)$$

### D. Constraints

*a) Generation Capacity Constraint*: For stable operation the real power output of each DG and ESS is restricted by lower and upper limits as follows:

$$P_{G_i}^{min} \leq P_{G_i} \leq P_{G_i}^{max}, \quad i=1,...N \quad (20)$$

$$SOC_{G_i}^{min} \leq SOC_{G_i} \leq SOC_{G_i}^{max}, \quad i=1,...M \quad (21)$$

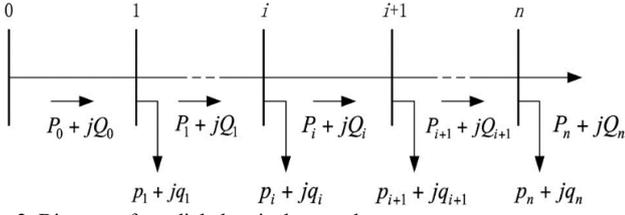

Fig. 2. Diagram of a radial electrical network.

*b) Power Balance Constraint*: The total electric power generation must cover the total electric power demand $P_D$ and the real power loss in distribution lines $P_{loss}$, hence:

$$\sum_{i=1}^{n} P_{G_i} - P_D - P_{loss} = 0 \qquad (22)$$

Consider an electrical network as shown in Fig. 2, calculation of implies solving the load flow problem, which has equality constraints on real and reactive power at each bus as follows:

$$P_{i+1} = P_i - r_i \left(P_i^2 + Q_i^2\right)/V_i^2 - P_{i+1} \qquad (23)$$

$$Q_{i+1} = Q_i - x_i \left(P_i^2 + Q_i^2\right)/V_i^2 - Q_{i+1} \qquad (24)$$

In the above equations, we assume $P_i$ is generated by both RES-based DG units which are subject to uncertainties and controllable DG units, $Q_i$ is generated by controllable DG units.

## V. Solution Algorithm

The cuckoo optimization algorithm (COA) is inspired by the life of a bird family, called Cuckoo [20]. Special lifestyle of these birds and their characteristics in egg laying and breeding has been the basic motivation for development of this new evolutionary optimization algorithm. Similar to other evolutionary methods, cuckoo optimization algorithm starts with an initial population. The cuckoo population, in different societies, is in two types: mature cuckoos and eggs. The effort to survive among cuckoos constitutes the basis of COA. During the survival competition some of the cuckoos or their eggs, demise. The survived cuckoo societies immigrate to a better environment and start reproducing and laying eggs. When moving toward goal point, the cuckoos do not fly all the way to the destination habitat. They only fly a part of the way and also have a deviation. This movement of cuckoos is clearly shown in Fig. 3 [20]. Cuckoos' survival effort hopefully converges to a state that there is only one cuckoo society all with the same profit values.

## VI. Results and Discussions

As shown in Fig. 4, the PG&E 69-bus distribution system with three type of DGs units and ESSs is used for solving the stochastic-based daily energy management (SDEM) problem in this paper. Details about the PG&E 69-bus test system can be found in [21].

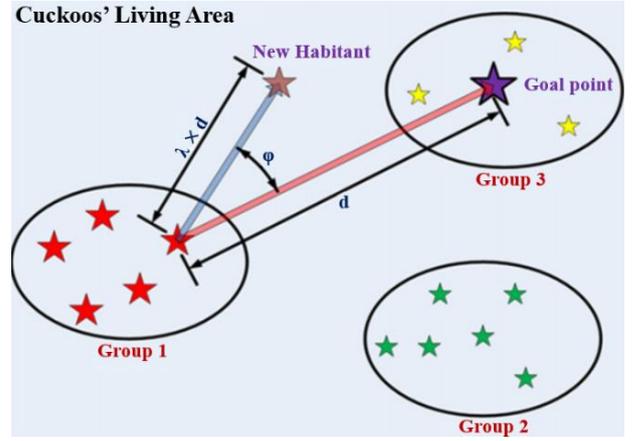

Fig. 3. Immigration of a sample cuckoo toward goal habitat.

In smart distribution network consists three type of SSERs and ESSs to power generation. All basic data for simulation are mentioned in the rest of this paper. The aim of this paper is to solve the optimal stochastic-based daily energy management of MGs at the presence of SSERs such as WT, PV, CHP, and energy storage system (ESS) for a given 24 hours. Because of the intermittent behavior of some renewable energy resources and load variation this paper is done based on uncertainty in input data. As an energy management this paper is solved the MG energy management problem based on measured data of 24h. In this paper, MG is modeled as a small disco with capability of power interchange between MG and the main grid. To consider the effect of power generation cost in MG, it is assumed that the characteristic of SSERs such as size, and technology of WT, PV, and CHP is different. In Table I, the values of parameters of WT and PV units are given that are used in simulation. Based on the above discussion, the cost of per kilowatt generated power by each SSERs at MG is different. Indeed, this paper is modeled as a daily energy management problem in market operation of smart distribution network. In this case, the stochastic energy management problem is solved using cuckoo optimization algorithm as a heuristic optimization algorithm.

In this paper, SDEM problem is solved considering PDF for some parameters of power generation. Proposed example solves SDEM problem at presence of load demand and renewable SSERs uncertainties. The results includes the total cost at any states of generated and transaction powers. The probabilistic analysis about SDEM has a good prospect in operation, planning, unit commitment problem studies, etc. PDF of the purchased and sold powers of MGs gives an insight vision to the dispatcher to evaluate the risk of change in system total cost with respect the variation in load and SSER power. This value at risk for purchase and sold power can be calculated and can be useful for dispatchers. The mean value of the PDF can use as the power with high probability of the appearance. Indeed, in this case, it is possible to evaluate the forecasting errors in load, wind, and solar short-term scheduling. The average value of generated power of PV, WT, CHP and ESS can be shown in Fig. 5 for 24 hours.

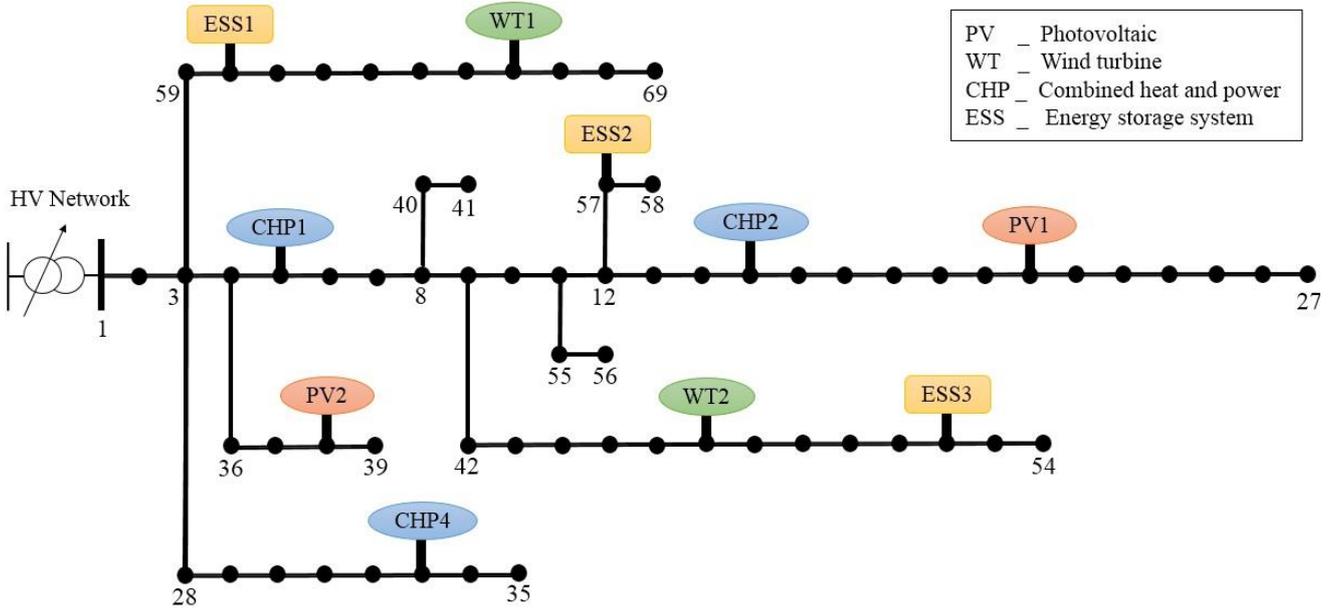

Fig. 4. Test distribution system with networked MGs

TABLE I. VALUES OF PARAMETERS IN WT AND PV.

| WT | $P_r$=250KW | $V_{CI}$=2m/s | $V_r$=14m/s | $V_{CO}$=25m/s |
|---|---|---|---|---|
| | n=4 | α=3 | β=12 | - |
| PV | $P_{STC}$=250KW | $G_{STC}$=1000W/m² | K=0.001 | $T_c$=25 C |

In probabilistic analysis PDF form is more accurate than mean value. PDF of ENS describes each possible state of it considering related probability. Based on this figure and Table II, the mean value of generated power of distributed generation including PV, WT and CHP with energy storage system at daily energy management in the distribution network. In this table each period is considered three hours. The voltage profile based on mean values has been shown in Fig. 6. According to this figure, with using the proposed method for daily energy management of MG, voltage profile of system is improved. Fig. 7 shows the power losses for proposed method and traditional energy management for smart distribution network. According to this figure, it is observed that by applying the proposed method for energy management the power losses of system are decreased. Based on correlation between network input and output variables, the output of the network parameters show probabilistic behavior. To explain this, for example, Fig. 8 is illustrated the probabilistic generated power of both WT and PV units in MG. In probabilistic analysis PDF form is more accurate than mean value. After obtaining the PDF of powers and reliability indices, the operation, emission and reliability costs can be represented.

In this paper one thousand scenarios are generated using Monte Carlo Simulation (MCS) to represent the prediction errors in the prediction horizon. Scenario reduction is applied to reduce the computation efforts while maintaining the solution accuracy. The 1000 generated scenarios are reduced to 30 scenarios.

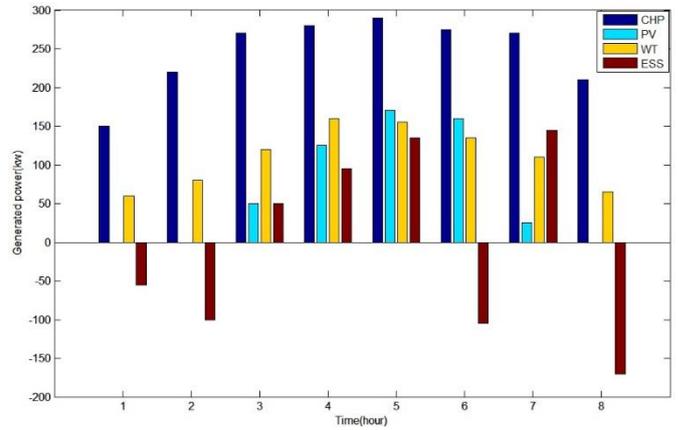

Fig. 5. Mean generated power of SSERs in DEM.

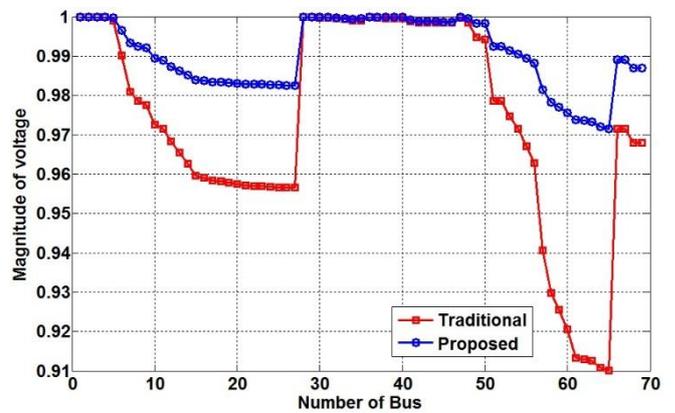

Fig. 6. Voltage profile of network based on mean value.

TABLE II. STATISTICAL ANALYSIS OF GENERATED POWERS.

| DER | PV | | WT | | CHP | | | ESS | | |
|---|---|---|---|---|---|---|---|---|---|---|
| Period (3h) | 1 | 2 | 1 | 2 | 1 | 2 | 3 | 1 | 2 | 3 |
| 1 | 0 | 0 | 26 | 34 | 52 | 43 | 55 | -15 | -23 | -18 |
| 2 | 0 | 0 | 33 | 47 | 73 | 68 | 80 | -27 | -39 | -35 |
| 3 | 31 | 19 | 60 | 60 | 78 | 95 | 93 | 22 | 15 | 13 |
| 4 | 68 | 55 | 85 | 78 | 85 | 96 | 98 | 27 | 38 | 35 |
| 5 | 89 | 76 | 66 | 89 | 88 | 102 | 100 | 34 | 42 | 58 |
| 6 | 84 | 78 | 53 | 59 | 94 | 95 | 87 | -29 | -33 | -37 |
| 7 | 12 | 11 | 44 | 67 | 79 | 104 | 92 | 54 | 45 | 42 |
| 8 | 0 | 0 | 23 | 38 | 66 | 75 | 65 | -58 | -63 | -49 |

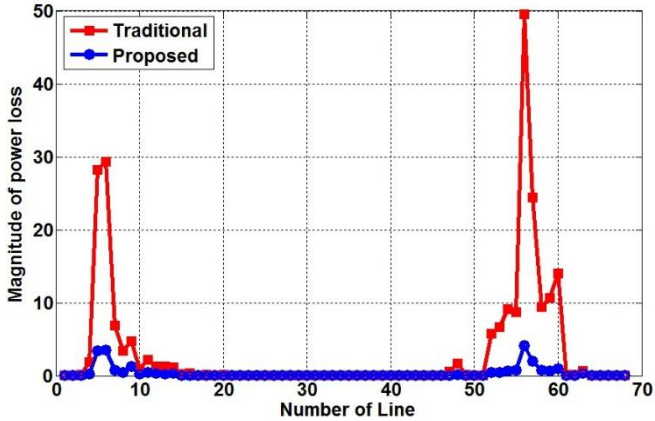

Fig. 7. Power losses of network based on mean value.

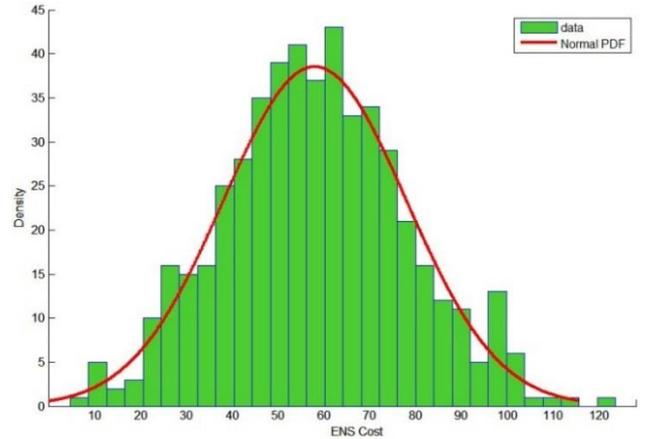

Fig. 9. PDF of ENS cost.

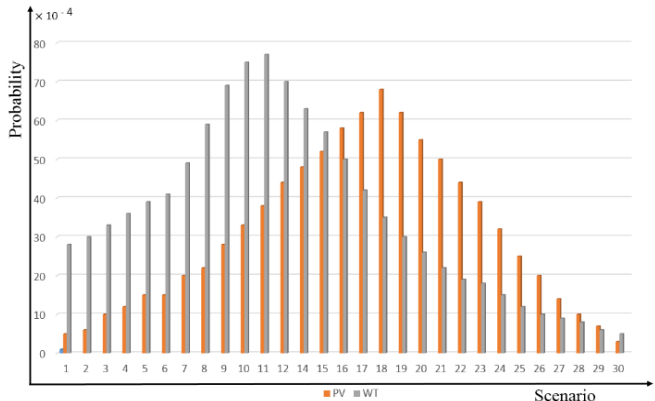

Fig.8. PDF of PV and WT generated power in MG.

The ENS cost of MG can be shown in PDF form in Fig. 9. In probabilistic analysis PDF form is more accurate than mean value. PDF of ENS describes each possible state of it considering related probability.

This paper presents a new stochastic energy scheduling scheme for the optimal energy management in a MG. In this scheme, the scheduling solution is computed based on the probabilistic forecasts of wind and solar power and MG load, thereby incorporating the uncertainties in both supply and demand.

The advantages for COA are less execution time, effective random numbers generation and organization because of having two absorbent operators and two efficient searching operators, the ability in selecting search area and in turn the ability for optimization of various problems, convergence to the absolutely identical solutions in each program iteration, not trapping in local minima and therefore high ability in extraction of global optimum points and needing no penalty factor in some problems. The results proved the robustness and effectiveness of the proposed algorithm and show that it could be used as a reliable tool for solving real and large scale EMS problems. Numerical simulation results show that, with the proposed scheme, the energy exchange between MG and utility grid is brought to a predefined trajectory. Comparison with some traditional scheduling schemes shows that, by incorporating the uncertainties in both sides of supply and demand, the proposed scheme exhibits significant improvement over the traditional schemes in terms of reference tracking.

## VII. Conclusion

In this paper a framework for stochastic optimal daily energy management of MGs in distribution systems considering economic, reliability and environment operation of distributed energy resources has been proposed. The uncertainty in MGs components such as DGs and load is modeled with PDF. The economic and reliable operation of MGs is formulated as an optimization problem. A stochastically and probabilistic modeling of both DGs and load is performed to determine the optimal management of MGs with minimum cost based on economic and reliability analysis of the power transaction between the MG and main grid. Based on the results, the mean, standard deviation, and PDF of each generated power with SSERs is determined considering optimization constraints. Statistical analysis for generated power and costs are given. The cuckoo optimization algorithm is applied to minimize the cost function as an optimization algorithm. Results show that it is possible to regulate the power demand and transaction between MGs and the main grid. Moreover, it is indicated that the power sharing between DGs with main grid can reduce the total operation and reliability cost of the future distribution network. One of the main results of this paper by probabilistic modeling of the input variables, the output variables can be represented as random variables. This leads to a better and comprehensive vision for network experts to manage the marginal operation of the network under uncertainties. The results indicate the ability of the proposed algorithm in finding the global optimum point of the functions for each run of the program. This can guarantee the robust energy management of smart distribution grids in the presence of network uncertainties. The obtained results show that the proposed framework can be considered as an efficient tool in optimal daily energy management of MGs in distribution networks with considering uncertainties.

There are two possible avenues for future work arising from this paper, namely: 1) Other uncertainty modeling approaches such as information gap decision theory, point estimate method, fuzzy mathematics, interval analysis, and robust optimization, and 2) considering other MGs management options like capacitor switching and network reconfiguration.